\documentclass[12pt]{article}

\def\bl{\begin{eqnarray}}
\def\el{\end{eqnarray}}
\def\bll{\begin{eqnarray*}}
\def\ell{\end{eqnarray*}}
\usepackage{amssymb}
\usepackage{amsmath,amsthm}
\usepackage{amsfonts}
\usepackage{fleqn}
\usepackage{titlesec}
\usepackage[backref=true,dvipdfm,draft=false,colorlinks=true,
bookmarksopen=true,bookmarks=true,bookmarksnumbered=false,
linkcolor=blue,citecolor=blue,urlcolor=blue,pdfauthor=YunBaoHuang]{hyperref}

\title{The subword complexity of  smooth words  on 2-letter alphabets}

\small \author{Yun Bao Huang \\
 Department of Mathematics\\
 Hangzhou Normal University\\
 Xiasha Economic Development Area\\
 Hangzhou, Zhejiang 310036, China\\
        huangyunbao@sina.com\\
        huangyunbao@gmail.com}
\date{2010.12.13}

\begin{document}
\numberwithin{equation}{section} \makeatletter
\titlelabel{\thetitle.\;}
\newtheorem{defn}{Definition}
\newtheorem{thm}[defn]{Theorem}
\newtheorem{lem}[defn]{Lemma}
\newtheorem{prop}[defn]{Proposition}
\newtheorem{cor}[defn]{Corollary}
\newtheorem{conj}[defn]{Conjecture}
\newtheorem{exmp}[defn]{Example}
\newtheorem{rem}[defn]{Remark}
\makeatother \maketitle

\begin{quote}
{\small {\bf Abstract.}  Let $\gamma_{a,b}(n)$ be the number of
smooth words of length $n$ over the alphabet $\{a,b\}$ with $a<b$.
Say that a smooth word $w$ is \emph{left fully extendable} (LFE) if
both $aw$ and $bw$ are smooth. In this paper, we prove that for any
positive number $\xi$ and positive integer $n_{0}$ such that the
proportion of $b$'s is larger than $\xi$ for each LFE word of length
exceeding $n_0$, there are two constants $c_{1}\,\textrm{and}\,
c_{2}$ such that for each positive integer $n$, one has
\begin{eqnarray}
c_{1}\cdot n^{\frac{\log (2b-1)}{\log
(1+(a+b-2)(1-\xi))}}<\gamma_{a,b}(n)< c_2\cdot n^{\frac{\log
(2b-1)}{\log (1+(a+b-2)\xi)}}.\nonumber
\end{eqnarray}
In particular, taking $a=1\text{ and }b=2$ in the above inequalities
arrives at Huang and Weakley's result. Moreover, for 2-letter even
alphabet $\{a,b\}$,  there are two suitable constants $c_1,\,c_2$
such that
\begin{eqnarray}
c_{1}\cdot n^{\frac{\log (2b-1)}{\log ((a+b)/2)}}<\gamma_{a,b}(n)<
c_2\cdot n^{\frac{\log (2b-1)}{\log ((a+b)/2)}}\textit{for each
positive integer $n$}.\nonumber
\end{eqnarray}
}
\end{quote}
\newpage
\section{Introduction}
The curious Kolakoski sequence $K$ which Kolakoski introduced in
~\cite{Kolakoski}, is the infinite sequence over the alphabet
$\Sigma=\{1,2\}$, which starts with 2 and equals the sequence
defined by its run lengths:
\begin{eqnarray}
K= \underbrace{22}_{2}\underbrace{11}_{2}\underbrace{2}_{1}
 \underbrace{1}_{1}\underbrace{22}_{2}\underbrace{1}_{1} \underbrace{22}_{2}\underbrace{11}_{2}
 \underbrace{2}_{1}\underbrace{11}_{2}\underbrace{22}_{2}\underbrace{\cdots}_{\cdots}.\nonumber
\end{eqnarray}
Here, a run is a maximal subsequence of consecutive identical
symbols. The Kolakoski sequence $K$ has received a remarkable
attention in $[1,2,3,\cdots,26]$. For research situations of the
Kolakoski sequence $K$ and related problems before 1996, readers can
refer to Dekking~\cite{Dekking3}.

Keane~\cite{Keane} asked whether the density of $1'$s in $K$ is 0.5.
Chv\'{a}tal~\cite{Chv} proved that the upper density of $1'$s as
well as the upper density of $2'$s in $K$ is less than 0.500838.
Steacy~\cite{Steacy} studied the structure in the kolakoski sequence
$K$ and obtained some conditions which are equivalent to Keane$'$s
problem.

In order to study wether the Kolakoski sequence $K$ is recurrent
and/or is closed under complement, Dekking~\cite{Dekking2}
introduced the notion of $C^\infty$-words over the alphabet
$\{1,2\}$ for the first time and noted that the finite factors of
$K$ must be $C^\infty$-words. Moreover, he proved that there exists
a suitable positive constant $c$ such that $c\cdot
n^{2.15}\leq\gamma(n)\leq n^{7.2}$ and conjectured that there are
suitable constants $c_1\text{ and }c_2$ such that $c_{1}n^{q}\leq
P_{K}(n)\leq c_{2}n^{q}$, where $\gamma(n)$ denote the number of
$C^{\infty}$-words of length $n$, $P_{K}(n)$ denote the number of
subwords (factors) of length $n$ which occur in the Kolakoski
sequence $K$,  $q=(\log 3)/\log (3/2)$.

Weakley~\cite{Weakly} showed that there are positive constants
$C_{1}$ and $C_{2}$ such that for each $n$ satisfying $B(k-1)+1\leq
n \leq A(k)+1$ for some $k$, $C_{1}n^{q}\leq \gamma(n) \leq
C_{2}n^{q}$, where $A(k),\; B(k)$ denote respectively  the minimum
and the maximal length of FE words of height $k$ (~\cite{Weakly}
Corollary 9).

Huang and Weakley~\cite{Huang4} proved that for any positive number
$\phi$ and positive integer $n_{0}$ satisfying
$|u|_{2}/|u|>\frac{1}{2}-\phi$ for each LDE word $u$ of length
exceeding $n_0$, there are two suitable constants
$c_{1}\,\textrm{and}\, c_{2}$ such that
\begin{eqnarray}
c_{1}n^{\frac{\log 3}{\log((3/2) + \phi + (2/N))}} < \gamma(n) <
c_{2}n^{\frac{\log 3}{\log((3/2) - \phi)}}\text{ for each }n\in N.
\nonumber
\end{eqnarray}
With the best value known for $\phi$, and large $N$, this gives
\begin{eqnarray}
c_{1}n^{2.7087} < \gamma_{1,2}(n) < c_{2}n^{2.7102}. \nonumber
\end{eqnarray}
A naturally arising question is whether or not we can establish the
estimates of subword complexity function of smooth words for the
other 2-letter alphabets. This paper is a study of subword
complexity function of smooth words for any 2-letter alphabets
(Theorem \ref{t1}). We establish the bounds of minimal and maximal
heights of smooth words of length $n$ (Lemma~\ref{Lem7}), the best
bounds of minimal and maximal heights of smooth words of length $n$
for 2-letter even alphabets (Lemma~\ref{Lem9}) and the good lower
and upper bounds of the subword complexity function
$\gamma_{a,b}(n)$ for 2-letter even alphabet $\{a,b\}$
(Theorem~\ref{t2}), which would give $\gamma_{a,b}(n)\approx c
n^{\log (2b-1)/\log \frac{a+b}{2}}$, where $c$ is a suitable
constant.

The paper is structured as follows. In Section \ref{s1}, we shall
first fix some notations and introduce some notions. Second in
Section \ref{s2}, we give some lemmas which are needed to establish
the estimates of the complexity function for arbitrary 2-letter
alphabets.  Third, in Section \ref{s3}, we obtain the lower and
upper bounds of the subword complexity function of smooth words.
Moreover, in Section \ref{s4}, we establish the good lower and upper
bounds of the subword complexity function $\gamma_{a,b}(n)$ for
2-letter even alphabets. Finally, in Section \ref{s5}, we end this
paper with some concluding remarks.
\section{Definitions and notation\label{s1}}
Let $\Sigma=\{a, b\}$ with $a<b$ and $a,\,b$ being positive
integers, $\Sigma^{*}$ denotes the free monoid over $\Sigma$ with
$\varepsilon$ as the empty word. A finite word over $\Sigma$ is an
element of $\Sigma^{*}$. If $w=w_{1}w_{2}\cdots$$w_{n}$, $w_{i}\in
\Sigma$ for $i=1, 2, \cdots, n$, then $n$ is called the length of
the word $w$ and is denoted by $|w|$. Let $|w|_{\alpha}$ be the
number of $\alpha$ which occur in $w$ for $\alpha\in \Sigma$, then
$|w|=|w|_{a}+|w|_{b}$.

  Given a word $w\in\Sigma^{*}$, a $\mathit{factor\;(\text{or}\; subword)}$ $u$ of $w$ is a word $u\in
\Sigma^{*}$ such that there exist $ x, y \in\Sigma^{*}$ such that
$w=xuy$.
   If $x=\varepsilon$ then $u$ is
called $\mathit{prefix}$. A $\mathit{run}$ (or $\mathit{block}$) is
a maximal factor of the form $u=\alpha^{k}, \alpha\in\Sigma$.
Finally, $N$ is the set of positive integers and the cardinal number
of $A$ is denoted by $|A|$ for a set $A$.

   The $\mathit{reversal}$ (or $\mathit{mirror\, image}$)\label{def4} of $u=u_{1}u_{2}\cdots$$u_{n}\in\Sigma^{*}$ is
the word $\tilde{u}=u_{n}u_{n-1}\cdots$$u_{2}\\u_{1}$. The
$\mathit{complement}$ (or $\mathit{permutation}$) of
$u=u_{1}u_{2}\cdots$$u_{n}\in\Sigma^{*}$ is the word
$\bar{u}=\bar{u}_{1}\bar{u}_{2}\cdots$$\bar{u}_{n}$, where
$\bar{a}=b, \bar{b}=a$.

    Now we generalize the definition of differentiable words, which Dekking first
introduced in ~\cite{Dekking2}, to over  arbitrary 2-letter alphabet
$\{a,\,b\}$ from the alphabet $\{1,\,2\}$.

To do so, for $w\in \Sigma^{*}$,  $r(w)$ denotes the number of runs
of $w$, $\mathit{fr}(w)$ and $\mathit{lr}(w)$ denote the first and
last runs of $w$ respectively, and $\mathit{lfr}(w)$ and
$\mathit{llr}(w)$ denote the lengths of the first and last run of
$w$ respectively. For example, if $w=a^{2}b^{2b}a^{a}b^{3}$, then
$r(w)=4$, $\mathit{fr}(w)=a^{2}$, $\mathit{lr}(w)=b^{3}$,
$\mathit{lfr}(w)=2$ and $\mathit{llr}(w)=3$.

Then we first need to introduce the concept of the closure of a word
$w$ over $\Sigma$ in order to establish the notion of differentiable
word for arbitrary 2-letter alphabets.
\begin{defn}\label{def1}
Let $w\in \Sigma^{*}$ and
\begin{eqnarray}
w&=&\alpha^{t_{1}}\bar{\alpha}^{t_{2}}\ldots
\beta^{t_{k}},\label{hk1}
\end{eqnarray}
$\textrm{where} \,\alpha\in \Sigma,\,\beta=\alpha\text{ if }2\nmid
k,\text{ or else } \beta=\bar{\alpha},\,1\leq t_{i}\leq b \,\,
\textrm{for} \,1\leq i\leq k.$
\begin{eqnarray}
\hat{w}=\left\{\begin{array}{ll}
w ,& \mathit{lfr}(w)\leq a \;\textrm{and}\;\mathit{llr}(w)\leq a\\
\alpha^{b-t_{1}}w,& \mathit{lfr}(w)> a \;\textrm{and}\;\mathit{llr}(w)\leq a\\
w\beta^{b-t_{k}}, & \mathit{lfr}(w)\leq a \;\textrm{and}\;\mathit{llr}(w)> a\\
\alpha^{b-t_{1}}w\beta^{b-t_{k}}, & \mathit{lfr}(w)> a \;\textrm{and}\;\mathit{llr}(w)> a\\
\end{array}\right.\nonumber
\end{eqnarray}
Then $\hat{w}$ is said to be the closure of a word $w$.
\end{defn}
For example, let $w=3311133313133311133,\, u=3313133311$, then $u$
is a factor of $w$, and $\hat{w}=333111333131333111333,\,
\hat{u}=333131333111$. Thus $\hat{u}$ is a factor of $\hat{w}$,
which also holds in general (see Lemma \ref{Lem2} (1)).
\begin{defn}\label{def2}
Let $w\in \Sigma^{*}$ be of the form  $(\ref{hk1})$. If the length
of every run of $w$ only  takes $a\textrm{ or } b$ except for the
lengths of the first and last runs, then we call that $w$ is
differentiable, and its derivative, denoted by $D(w)$, is the word
whose $j$th symbol equals the length of the $j$th run of $w$,
discarding the first and/or the last run if its length is less than
$b$.
\end{defn}
If $\hat{w}$ is differentiable, then we call that $w$ is closurely
differentiable. If a finite word $w$ is arbitrarily often closurely
differentiable, then we call $w$ a $C^{\infty}_{a,b}$-word or a
smooth word over the alphabet $\{a,b\}$, and the set of all smooth
words over the alphabet $\{a,b\}$ is denoted  by $C^{\infty}_{a,b}$
or $C^{\infty}$.

Let $\rho(w)=D(\hat{w})$\label{def6}, then it is clear that $w$ is a
smooth word if and only if there is a positive integer $k$ such that
$\rho^{k}(w)=\varepsilon$.

Note that if $b=a+1$ then $\hat{w}=w$. Thus, $w$ is differentiable
if and only if $w$ is closurely differentiable, which suggests that
$w$ is a smooth word if and only if there is a positive integer $k$
such that $D^{k}(w)=\varepsilon$.

By the definition 2, it is clear that if $b-a\geq 2$ and $a\neq 2$,
then $a^{b-1}b^{a}a^{a}b^{b-1}$ is differentiable but not closurely
differentiable. Moreover, D is an operator from $\Sigma^{*}$ to
$\Sigma^{*}$, $r(w)\leq |D(w)|+2$\label{huang} and\\
\begin{eqnarray}\label{hk2}
D(\hat{w})=\left\{\begin{array}{ll}
bD(w),& b>\mathit{lfr}(w)> a \;\textrm{and}\;\mathit{llr}(w)\leq a\\
D(w)b, &  b>\mathit{llr}(w)> a \,\;\textrm{and}\;\mathit{lfr}(w)\leq a\\
bD(w)b, &  b>\mathit{lfr}(w)> a \;\textrm{and}\;b>\mathit{llr}(w)> a\\
D(w) ,& \textrm{otherwise}
\end{array}\right..
\end{eqnarray}
From  (\ref{hk2}), it follows that if $w$ is closurely
differentiable, then it must be differentiable.

A word $v$ such that $D(v)=w$ is said to be a primitive of $w$. The
two primitives of $w$ having minimal length are the shortest
primitives of $w$. For example, $b$ have $2b^{2}$ primitives of the
form $\alpha^{i}\bar{\alpha}^{b}\alpha^{j}$, where $\alpha=a,b$, $i,
j=0,1,\cdots b-1$, and $a^{b},\,b^{b}$ are the shortest primitives.
It is easy to see that for any word $w \in C^{\infty}$, there are at
most $2b^{2}$ primitives, and the difference of lengths of two
primitives of $w$ is at most $2(b-1)$.

The height of a smooth word $w$ is the smallest integer $k$ such
that $D^{k+1}(w)=\varepsilon$. We write $ht(w)$ for the height of
$w$. For example, if
$w=32^{3}3^{3}2^{3}3^{2}2^{2}3^{2}2^{3}3^{3}2^{3}3$, then $ht(w)=3$.

It immediately follows from the definition \ref{def2} that

(1) $D(\tilde{u})=\widetilde{D(u)},\; D(\bar{u})=D(u)$ for each
$u\in\Sigma^{*}$

(2) $w \in C^{\infty}\Longleftrightarrow\bar{w}, \tilde{w}\in
C^{\infty}$.
\section{Some lemmas\label{s2}}
The following Lemmas \ref{Lem2} to \ref{Lem4} reveal the relations
among the operators mirror image, complement, closure and
derivative.
\begin{lem}[\cite{Huang6}, Lemma 5]\label{Lem2}
Let $w$ be a differentiable word  and $u$ is a factor of $w$. Then

\rm (1) \it both  $\hat{u}$ and $w$ are factors of $\hat{w}$;

\rm (2) \it $\hat{\tilde{w}}=\tilde{\hat{w}},
\hat{\bar{w}}=\bar{\hat{w}}$;

\rm (3) \it $D(u)$ is a factor of $D(w)$;

\rm (4) \it If $w$ is closurely differentiable, then  both $\rho(u)$
and $D(w)$ are factors of $\rho(w)$, and $\rho(\bar{w})=\rho(w)$,
$\rho(\tilde{w})=\widetilde{\rho(w)}$.
\end{lem}
\noindent{\bf Proof. } {(1) From the definition \ref{def1} of the
closure of a word, it follows the assertion (1).

(2) It immediately follows from the definitions of the closure,
complement and mirror image of a word $w$ and the definition of the
operators $\rho$.

(3) Since $u$ is a factor of $w$, by the definition \ref{def2} of
the derivative of a word $w$, we see that $D(u)$ is a factor of
$D(w)$.

(4) Since $w$ is closurely differentiable and $\rho(w)=D(\hat{w})$,
by the assertion (1), $\hat{u}$ and $w$ are both factors of
$\hat{w}$. Moreover by the assertion (3), we see that $D(\hat{u})$
and $D(w)$ are factors of $D(\hat{w})$, that is, both $\rho(u)$ and
$D(w)$ are factors of $\rho(w)$. Finally, by the assertion (2), we
have
$\rho(\bar{w})=D(\hat{\bar{w}})=D(\bar{\hat{w}})=D(\hat{w})=\rho(w)$.
Similarly,
$\rho(\tilde{w})=D(\hat{\tilde{w}})=D(\tilde{\hat{w}})=\widetilde{D(\hat{w})}=\widetilde{\rho(w)}$.
$\Box$ }

 From the definitions 1-2, it immediately follows that
\begin{lem}[\cite{Huang6}, Lemma 6]\label{Lem3}
Let $w=w_{1}w_{2}\cdots w_{n}$ be a differentiable word with
$n\geq a+1$.

\rm (1) \it If \,$\mathit{lfr}(w)=b$ then $w_{1}w$ is not a
differentiable word and $D(\bar{w}^{i}_{1}w)=D(w)$ for $i\leq b-1$;

\rm (2) \it If \,$\mathit{lfr}(w)<b$ then
$D(w^{b-\mathit{lfr}(w)}_{1}w)=bD(w)$;

\rm (3) \it If \,$\mathit{lfr}(w)\leq a$ and $r(w)>1$ then
$D(\bar{w}_{1}w^{a-\mathit{lfr}(w)}_{1}w)=aD(w)$.\;\;$\Box$
\end{lem}
\begin{lem}[\cite{Huang6}, Lemma 7]\label{Lem4}\rm
(1) \it Let $w=w_{1}w_{2}\cdots w_{n}$ be a
 smooth word. Then any factor of $w$ is also a
 smooth  word;\rm

 (2) \it Any smooth  word $w=w_{1}w_{2}\cdots w_{n}$ has both a left and a right
smooth extensions.
\end{lem}
\noindent{\bf Proof. } (1) If $w$ is a smooth word and $u$ is a
factor of $w$, then note that $w\in C^{\infty}\Longleftrightarrow
\rho^{k}(w)=\varepsilon$ for some positive integer $k$, by Lemma
\ref{Lem2} (4), we obtain that $\rho^{i}(u)$ is a factor of
$\rho^{i}(w)$ for any positive integer $i\leq k$. And hence
$\rho^{k}(w)=\varepsilon$ suggests $\rho^{k}(u)=\varepsilon$, so
that $u$ is a smooth  word.

(2) We verify the assertion (2) by induction on $|w|$. Since
$D(\tilde{w})=\widetilde{D(w)}$, we only need to verify that $w$ has
a left smooth  extension. It is clear that if $r(w)\leq 1$, where
$r(w)$ is the number of runs of $w$, then the assertion (2) holds.
We proceed to the induction step. Assume now that $r(w)\geq 2$ and
the assertion (2) holds for smooth  words shorter than $w$.

If $\mathit{lfr}(w)\leq a$ then by Lemma \ref{Lem3} (2-3), we have
$D(\bar{w}_{1}w^{a-\mathit{lfr}(w)}_{1}w)=aD(w)$ and
$D(w^{b-\mathit{lfr}(w)}_{1}w)=bD(w)$. Thus by $|D(w)|<|w|$, we see
that at least one of $aD(w)$ and $bD(w)$ is a smooth  word, which
means that $w$ has a left smooth  extension.

If $b>\mathit{lfr}(w)> a$, then by $w\in C^{\infty}$, we obtain that
$\hat{w}$ is a left smooth  extension of $w$.

If $b=\mathit{lfr}(w)$, then by Lemma \ref{Lem3} (1), we see that
$\bar{w}_{1}w$ is a left smooth  extension of $w$. $\Box$

Now we are in a position to generalize the notion of LDE words to
over arbitrary 2-letter alphabets from the alphabet $\{1,2\}$, which
Weakley first introduced in \cite{Weakly}.

If $aw\text{ and } bw$ are  both smooth, then the word $w$ is said
to be \emph{left fully extendable} (LFE). Clearly, LFE words are
closed under complement. For every nonnegative integer $k$, let
$LF_k$ denote the set of LFE words of length $k$.

Let $\gamma_{a,b}(k)$ denote the number of smooth words of length
$k$ over the alphabet $\{a,b\}$. Being similar to Weakley
\cite{Weakly}, define the differences of $\gamma_{a,b}$ by
$\gamma^{'}_{a,b}(k)=\gamma_{a,b}(k+1)-\gamma_{a,b}(k)$ for each
$k\geq 0$. From the definition of LFE words, it immediately follows
that $\gamma^{'}_{a,b}(k)=|LF_k|$ for each nonnegative integer $k$.
Since $\gamma(0)=\gamma^{'}(0)=1$, so
\begin{eqnarray}\label{eq1}
\gamma_{a,b}(k)=\gamma_{a,b}(0)+\sum^{k-1}_{i=0}\gamma^{'}_{a,b}(i)=1+\sum^{k-1}_{i=0}|LF_i|=2+\sum^{k-1}_{i=1}|LF_i|
\text{ for }k\geq 1.
\end{eqnarray}
\begin{lem}\label{Lem5}
\it Let $w=w_{1}w_{2}\cdots w_{k}$ be a
 smooth word, where $k\in N$. If $w$ is a LFE word then $D(w)$ is also  a
 LFE word, and if $k\geq b$ or $r(w)>1$ then $w=w_1^aw_{a+1}\ldots w_{k}$, where $w_1\neq
 w_{a+1}$.
\end{lem}
\noindent{\bf Proof.} Assume that $w$ is a LFE word of length
exceeding 0. If $k=|w|<b$ then it follows from both $w_1w$ and
$\bar{w}_1w$ being smooth words that $w=\overbrace{\alpha^a\ldots
\beta^a}^t\bar{\beta}^j$, where $\alpha\in \Sigma$, $\beta=\alpha$
if $2\nmid t$, otherwise $\beta=\bar{\alpha}$, $0\leq j,\,t\leq
b-1,\;j+ t\geq 1,\;k=t\cdot a+j$. So $D(w)=a^{t-1}$ if $j,\,t\geq
1$, or else $D(w)=\varepsilon$. So, in view of $t<b$ we see that
$aD(w)\text{ and }bD(w)$ are both smooth words, that is, $D(w)$ is a
LFE word.

If $k\geq b$, since $w_1w$ is a smooth word, we get
$\mathit{lfr}(w)< b$, which suggests that $w=w_1^aw_{a+1}\ldots
w_{k}$ and $w_{a+1}\neq w_1$ by $\bar{w}_1w\in C^\infty$. Moreover,
note that each smooth word has a left smooth extension
(Lemma~\ref{Lem4} (2)), from $w_1w,\bar{w}_1w\in C^\infty$ it
follows that $aD(w)(=D(\bar{w}_1w))\text{ and
}bD(w)(=D(w_1^{b-a}w))$ are both smooth words, that is, $D(w)$ is a
LFE word. $\Box$

Let $LF$ denote the set $\bigcup\limits^{\infty}_{i=0}LF_i$ and
$P(A)=\{u\in LF: |u|>0\text{ and }D(u)\in A\}$ for $A\subseteq
\Sigma^*$. We now give the number of the elements contained in
$P^j(\varepsilon)$ for $j\in N$.
\begin{lem}\label{Lem6}
$|P^j(\varepsilon)|=4(b-1)(2b-1)^{j-1} \text{ for } j\in N.$
\end{lem}
\noindent{\bf Proof.} By Lemma \ref{Lem5} and the definition of
$P(A)$, we see that $P^{j+1}(\varepsilon)$ is exactly composed of
all LFE primitives of $P^j(\varepsilon)$.

Since for each LFE words of the form $\alpha\ldots b$ there are
exactly $2b$ LFE primitives:
\begin{eqnarray}
\bar{\beta}^a\Delta_\beta^{-1}(\alpha\ldots b)\gamma^j, \nonumber
\end{eqnarray}
where $ \alpha,\,\beta\in \Sigma,\;j=0,1,\ldots,b-1;
\,\gamma=\beta\text{ if }2\mid |\alpha\ldots b|, \text{ or else
}\gamma=\bar{\beta}$.

for each LFE words of the form $\alpha\ldots a$ there are exactly
$2(b-1)$ LFE primitives:
\begin{eqnarray}
\bar{\beta}^a\Delta_\beta^{-1}(\alpha\ldots a)\gamma^j,\nonumber
\end{eqnarray}
where $ \alpha,\,\beta\in \Sigma,\;j=1,\ldots,b-1;
\,\gamma=\beta\text{ if }2\mid |\alpha\ldots a|, \text{ or else
}\gamma=\bar{\beta}$.

In addition, because of $\overline{\alpha\ldots
b}=\bar{\alpha}\ldots a$, we see that the numbers of LFE words of
the form both $\alpha\ldots b$ and $\alpha\ldots a$ are equal in all
LFE words of the same heights. It follows that
\begin{eqnarray}
|P^j(\varepsilon)|&=&2b\cdot
\frac{1}{2}|P^{j-1}(\varepsilon)|+2(b-1)\cdot \frac{1}{2}
|P^{j-1}(\varepsilon)|\nonumber\\
&=&(2b-1)|P^{j-1}(\varepsilon)|\text{ for }j\in N,\nonumber
\end{eqnarray}
which suggests that
\begin{eqnarray}\label{eqq2}
|P^j(\varepsilon)|=(2b-1)^{j-1}|P(\varepsilon)|.
\end{eqnarray}
Since the primitives of $\varepsilon$ are of the form
$\alpha^i\bar{\alpha}^j$, where $0\leq i,j\leq b-1$ and  $i+j\geq
1$, so  by $\bar{\alpha}(\alpha^i\bar{\alpha}^j)\in C^\infty$, we
get that if $i\geq 1\text{ and }j\geq 1$ then $i=a,\,j=1,2,\ldots,
b-1$, if $i\geq 1\text{ and }j=0$ then
$\alpha^{i+1},\,\bar{\alpha}\alpha^i\in C^\infty$, which suggests
$1\leq i\leq b-1$. Thus $\varepsilon$ have exactly $4(b-1)$ LFE
primitives. Thus (\ref{eqq2}) gives the desired result. $\Box$
\begin{lem} \label{l3} Let $\xi$ be a positive real number
 and $n_{0}$ a positive integer such that
\begin{eqnarray}\label{eqq3}
|u|_b/|u| > \xi \text{ for every LFE word $u$ of length exceeding
 }n_{0}.
\end{eqnarray}
Then

$(1)$  $|D(w)|\leq \alpha |w|$ for each LFE word $w$ with
$|w|>N_{0}$, where $\alpha=1/(1+(a+b-2)\xi),\,N_{0}$ is a suitable
positive integer.

$(2)$ $|w|\leq \beta |D(w)|+q$ for each LFE word $w$, where
$\beta=1+(a+b-2)(1-\xi),\,q$ is a suitable positive constant.
\end{lem}
\noindent{\bf Proof.} (1) Since the complement of any smooth word is
still a smooth word of the same length and $|u|_a=|\bar{u}|_b$, the
hypothesis (\ref{eqq3}) of Lemma~\ref{l3} means that
\begin{eqnarray}\label{eqq1}
|u|_a/|u| > \xi \text{ for every LFE word $u$ with } |u| \geq n_{0}.
\end{eqnarray}

It is easy to see
\begin{equation}\label{eq4}
|w|=|D(w)|+(a-1)|D(w)|_a+(b-1)|D(w)|_b+c,\text{ where }0\leq c\leq
2(b-1).
\end{equation}
From (\ref{eqq3}) to (\ref{eq4}), one has $|w|\geq
(1+(a+b-2)\xi)|D(w)| \text{ for } |D(w)|_b/|D(w)|>\xi$, which
implies $|D(w)|<\alpha |w|$ for every LFE word $w$ with $|w|\geq
N_{0}$, where $N_{0}$ is a suitable positive integer such that
$|D(w)|\geq n_0$ as soon as $|w|\geq N_{0}$.

(2)    As $|D(w)|_a/|D(w)|+|D(w)|_b/|D(w)|=1$, from (\ref{eqq3}) and
(\ref{eqq1}) ones get
\begin{eqnarray}
|D(w)|_a/|D(w)|& < & 1 - \xi\text{ for each LFE word $w$ with } |w| \geq N_{0},\label{eqq4}\\
|D(w)|_b/|D(w)|& < & 1 - \xi\text{ for each LFE word $w$ with } |w|
\geq N_{0}.\label{eqq5}
\end{eqnarray}
So, from (\ref{eq4}) to (\ref{eqq5}) it follows that $|w|\leq \beta
|D(w)|+2(b-1)$ for $|w| \geq N_{0}$, which means that (2) also
holds. $\Box$

The next lemma establishes the bounds of the heights of
$C^\infty$-words of length $n$, which is of independent interest.
\begin{lem}\label{Lem7}  Let $ht_{max}(n)$ and $ht_{min}(n)$ denote respectively the maximal
and the minimal heights of LFE words of length $n$, then for any
positive number $\xi$ and positive integer $n_{0}$ satisfying
$|u|_b/|u|>\xi$ for each LFE word $u$ with $|u|>n_{0}$, there are
two suitable constants $t_1$ and $t_2$ such that for every positive
integer $n$, one has
\begin{eqnarray}
ht_{min}(n)&>& \frac{\log n}{\log (1+(a+b-2)(1-\xi))}+t_1
\label{eq2},\\
ht_{max}(n)&<& \frac{\log n}{\log (1+(a+b-2)\xi)}+t_2,\label{eq3}
\end{eqnarray}
where $t_1=-\frac{\log (2(b-1)+\frac{q}{\beta-1})}{\log \beta}$, $q$
and $\beta$ are determined by Lemma~\ref{l3} $(2)$.
\end{lem}
\noindent{\bf Proof.} First, one checks (\ref{eq3}). Since $|D(w)|<
|w|$ for each $|w|>0$, and $|D(w)|\leq \alpha |w|$ for each LFE word
$w$ satisfying $|w|\geq N_{0}$ by Lemma \ref{l3} (1).

Let $k_{0}-1$ be the greatest height of all LFE words of
length$<N_{0}$ and $m_0$ is the least positive integer such that if
$|w|=m_0$, then the height of each LFE  word $w$ is no less than
$k_{0}$. Thus for every LFE word $w$, if $|w|\geq m_0$, then one can
get
\begin{eqnarray}
|D^{k}(w)|< \alpha^{k-k_{0}}|w| \text{ for }k\geq k_{0}.\nonumber
\end{eqnarray}
Hence
\begin{eqnarray}
|D^{k}(w)|<1&\text{as soon as}&\alpha^{k-k_{0}}|w|\leq
1\nonumber\\
&\Longleftrightarrow&k\geq \log(|w|)/\log (1/\alpha)+k_0,\nonumber
\end{eqnarray}
which means that the height $k-1$ of $w$ is smaller than $\log
(|w|)/\log (1+(a+b-2)\xi)+k_0$. Since there are only finite many LFE
words satisfying $|w|< m_0$, so there is a suitable constant $t_2$
such that (\ref{eq3}) holds for each LFE word.

Second, by Lemma \ref{l3} (2), one has $|w|\leq \beta |D(w)|+q$ for
each LFE word $w$, where $\beta=1+(a+b-2)(1-\xi),\,q$ is a suitable
constant, which means that
\begin{eqnarray}
|w|&<&\beta^{k} |D^{k}(w)|+q\frac{\beta^{k}-1}{\beta-1}\nonumber\\
&<&2(b-1)\beta^{k}+\frac{q\beta^{k}}{\beta-1}\nonumber\\
&=&(2(b-1)+\frac{q}{\beta-1})\beta^{k}\nonumber\\
&=&m\beta^{k},\nonumber
\end{eqnarray}
where $m=2(b-1)+q/(\beta-1)$, $k$ is the height of $w$. Thus the
length $|w|$ of a LFE word $w$ of height $k$ is less than
$m\beta^{k}$, and it follows that
\begin{eqnarray}
k>(\log |w|-\log m)/\log \beta, \nonumber
\end{eqnarray}
which gives the
desired lower bound of $ht_{min}(n)$, where $t_1=-\log m/\log \beta$. $\Box$\\
{\bf Remark 1.} (1) From (\ref{eqq3}) and (\ref{eqq1}) it
immediately follows that the positive real number $\xi$ satisfying
the condition (\ref{eqq3}) must be smaller than 1/2.

(2) From the proof of Lemma~\ref{l3} we easily see that if we
substitute  LFE words in Lemma~\ref{l3} with some infinite subclass
of smooth words, which is closed under complement, then the
corresponding result also holds.

(3) From the proof of Lemma~\ref{Lem7} we see that if we replace
 LFE words in Lemma \ref{Lem7} with some infinite subclass of
smooth words, which is closed under both complement and the operator
$D$, then the corresponding result still holds.

\section{The subword complexity of smooth words\label{s3}}
Now, we can establish our main result on subword complexity function
$\gamma_{a,b}(n)$ of smooth words over 2-letter alphabets.
\begin{thm}\label{t1}
For any positive real number $\xi$ and positive integer $n_{0}$
satisfying $|u|_b/|u|>\xi$ for every LFE word $u$ with $|u|>n_0$,
there exist two suitable constants $c_{1}\text{ and }\, c_{2}$ such
that
\begin{eqnarray}
c_{1} n^{\frac{\log (2b-1)}{\log
(1+(a+b-2)(1-\xi))}}\leq\gamma_{a,b}(n)\leq c_{2} n^{\frac{\log
(2b-1)}{\log (1+(a+b-2)\xi)}}\nonumber
\end{eqnarray}
for every positive integer $n$.
\end{thm}
\noindent {\bf Proof.} First, from the definition of $ht_{max}(n)$,
one sees that the length of LFE words of the height larger than
$ht_{max}(n)$ must be larger than $n$. Thus $\cup_{i=1}^{n-1} LF_{i}
\subseteq \cup_{j=1}^{ht_{max}(n)} P^{j}(\varepsilon)$. So from
(\ref{eq1}) and Lemma \ref{Lem6}, for any $n\in N$, one has
\begin{eqnarray}
\gamma_{a,b}(n) &=& 2 + \sum_{i=1}^{n-1} |LF_{i}|\nonumber\\
 &\leq& 2 + \sum_{j=1}^{ht_{max}(n)} |P^{j}(\varepsilon)|\nonumber\\
 & = &2 + \sum_{j=1}^{ht_{max}(n)} 4(b-1) \cdot (2b-1)^{j-1}\nonumber\\
& = & 2 \cdot (2b-1)^{ht_{max}(n)}.\label{eeq1}
\end{eqnarray}
So combining (\ref{eq3}) and (\ref{eeq1}) yields the desired upper
bound of $\gamma_{a,b}(n)$, where $c_2=2(2b-1)^{t_2}$.

Second, from the definition of $ht_{min}(n)$, it follows that the
length of all LFE words with the height no more than $ht_{min}(n)-1$
must be less than $n$. Thus, again from (\ref{eq1}) and Lemma
\ref{Lem6}, for any $n\in N$ one can get
\begin{eqnarray}
\gamma_{a,b}(n) &=& 2 + \sum_{i=1}^{n-1} |LF_{i}|\nonumber\\
 &\geq& 2 + \sum_{j=1}^{k} |P^{j}(\varepsilon)|\nonumber\\
 & = &2 + \sum_{j=1}^{k} 4(b-1) \cdot (2b-1)^{j-1}\nonumber\\
& = & 2 \cdot (2b-1)^{k},\label{eeq2}
\end{eqnarray}
where $k=ht_{min}(n)-1$. Thus, the desired lower bound of
$\gamma_{a,b}(n)$ is obtained from (\ref{eq2}) and (\ref{eeq2}),
where $c_1=2(2b-1)^{t_1-1}$, $t_1$ is decided by Lemma~\ref{Lem7}. $\Box$\\
{\bf Remark 2.} Theorem~\ref{t1} indicates that only if we could get
lower and upper bounds of letters frequency of LFE words, then
correspondingly we could obtain an estimate of subword complexity
function $\gamma_{a,b}(n)$ of smooth words.

Taking $\Sigma=\{1,2\}$ in Theorem \ref{t1}, we obtain
\begin{cor}\label{c1}
For any positive number $\xi$ and positive integer $n_{0}$
satisfying $|u|_2/|u|>\xi$ for each LDE word $u$ with $|u|>n_0$,
there exist two suitable constants $c_{1}\text{ and }\,c_{2}$ such
that
\begin{eqnarray}
c_{1}\cdot n^{\frac{\log 3}{\log (2-\xi)}}\leq\gamma_{1,2}(n)\leq
c_{2}\cdot n^{\frac{\log 3}{\log (1+\xi)}} \text{ for each } n\in
N.\nonumber
\end{eqnarray}
\end{cor}
It is obvious that Corollary~\ref{c1} suggests the main Theorem 1 in
\cite{Huang4}.
\section{The subword complexity of smooth words on 2-letter even alphabets\label{s4}}
\begin{lem} \label{Lem8} If $w$ is a 2-times differentiable finite word over 2-letter even alphabet
$\{a,b\}$, then
\begin{eqnarray}
(1)\;||w|_{a}-|w|_{b}|\leq b;\nonumber
\end{eqnarray}
\begin{eqnarray}
(2)\;\frac{1}{2}-\frac{b}{2|w|}\leq \frac{|w|_b}{|w|}\leq
\frac{1}{2}+\frac{b}{2|w|}; \nonumber
\end{eqnarray}
\begin{eqnarray}
(3)\;\lim_{|w|\rightarrow \infty}
\frac{|w|_{a}}{|w|}=\lim_{|w|\rightarrow \infty} \frac{|w|_{b}}{|w|}
                             =\frac{1}{2};
\nonumber
\end{eqnarray}
\begin{eqnarray}
(4)\;\rho |D(w)|-q_2\leq |w| \leq \rho |D(w)|+q_1, \nonumber
\end{eqnarray}
where $q_1=(\rho-1)b+2(b-1)$, $q_2=(\rho-1)b$, $\rho=\frac{a+b}{2}$.
\end{lem}
\noindent{\bf Proof.} It is obvious that $(1)\Rightarrow
(2)\Rightarrow (3)$. So we only need to check (1) and (4).

(1) Since $w\in C^{2}_{a,b}$, we have $D^{2}(w)\in \Sigma^{*}$. Thus
$D^{2}(w)=\alpha^{t_{1}}\bar{\alpha}^{t_{2}}\cdots \beta^{t_{k}}$,
where $\alpha\in \Sigma$, $t_{i}\in N$ for $i=1,2,\cdots,k$, and if
$2\mid k$ then $\beta=\bar{\alpha}$, otherwise $\beta=\alpha$. It
follows that
\begin{eqnarray}
\Delta^{-1}_{\gamma_{1}}(D^{2}(w))&=&\overbrace{\gamma^{\alpha}_{1}\bar{\gamma_{1}}^{\alpha}\cdots\bar{\gamma_{2}}^{\alpha}}^{t_{1}}
\overbrace{\gamma_{2}^{\bar{\alpha}}\bar{\gamma_{2}}^{\bar{\alpha}}\cdots\bar{\gamma_{3}}^{\bar{\alpha}}}^{t_{2}}\cdots
\overbrace{\gamma_{k}^{\beta}\bar{\gamma_{k}}^{\beta}\cdots\bar{\gamma}_{k+1}^{\beta}}^{t_{k}};\label{eq6.1}\\
D(w)&=&\bar{\gamma}^{i}_{1}\Delta^{-1}_{\gamma_{1}}(D^{2}(w))\gamma^{j}_{k+1}\label{eq6.2}
\end{eqnarray}
where $0\leq i,\,j\leq b-1$, $\gamma_{i}\in \Sigma$ and if $2\mid
t_{m}$ then $\gamma_{m+1}=\gamma_{m}$ or else
$\gamma_{m+1}=\bar{\gamma}_{m}$.

Note that $a$ and $b$ are both even numbers, from (\ref{eq6.1}) it
immediately follows
\begin{eqnarray}
\Delta^{-1}_{\gamma}(\Delta^{-1}_{\gamma_{1}}(D^{2}(w)))=\overbrace{\gamma^{\gamma_1}\bar{\gamma}^{\gamma_1}\cdots\bar{\gamma}^{\gamma_1}}^{\alpha}
\overbrace{\gamma^{\bar{\gamma_1}}\bar{\gamma}^{\bar{\gamma_1}}\cdots\bar{\gamma}^{\bar{\gamma_1}}}^{\alpha}\cdots
\overbrace{\gamma^{\bar{\gamma}_{k+1}}\bar{\gamma}^{\bar{\gamma}_{k+1}}\cdots\bar{\gamma}^{\bar{\gamma}_{k+1}}}^{\beta},\label{eq6.3}
\end{eqnarray}
where $\alpha,\beta,\gamma,\gamma_1,\cdots,\gamma_{k+1}\in \Sigma$.

Then (\ref{eq6.3}) gives
\begin{eqnarray}
|\Delta^{-1}_{\gamma}(\Delta^{-1}_{\gamma_{1}}(D^{2}(w)))|_a=|\Delta^{-1}_{\gamma}(\Delta^{-1}_{\gamma_{1}}(D^{2}(w)))|_b=\frac{1}{2}\label{eq6.4}.
\end{eqnarray}
Now from (\ref{eq6.2}) ones get
\begin{eqnarray}
w&=&\xi^{c_1}\Delta_{\bar{\xi}}^{-1}(D(w))\eta^{c_2}\nonumber\\
&=&\xi^{c_1}\Delta_{\bar{\xi}}^{-1}(\bar{\gamma}^i_1)\Delta_{\mu}^{-1}(\Delta^{-1}_{\gamma_{1}}(D^{2}(w)))\Delta_{\mu}^{-1}(\gamma_{k+1}^j)\eta^{c_2}\label{eq6.5}
\end{eqnarray}
where $0\leq i,j,c_1,c_2\leq b-1$, $\mu=\bar{\xi}$ if $2\mid i$ or
else $\mu=\xi$, $\eta=\bar{\xi}$ if $2\mid (i+j)$ or else
$\eta=\xi$. Note that
\begin{eqnarray}
||\xi^{c_1}\Delta_{\bar{\xi}}^{-1}(\bar{\gamma}^i_1)|_a-|\xi^{c_1}\Delta_{\bar{\xi}}^{-1}(\bar{\gamma}^i_1)|_b|\leq
b,\nonumber
\end{eqnarray}
\begin{eqnarray}
||\Delta_{\mu}^{-1}(\gamma_{k+1}^j)\eta^{c_2}|_a-|\Delta_{\mu}^{-1}(\gamma_{k+1}^j)\eta^{c_2}|_b|\leq
b.\nonumber
\end{eqnarray}
And if
\begin{eqnarray}
|\xi^{c_1}\Delta_{\bar{\xi}}^{-1}(\bar{\gamma}^i_1)|_\alpha\geq|\xi^{c_1}\Delta_{\bar{\xi}}^{-1}(\bar{\gamma}^i_1)|_{\bar{\alpha}}\nonumber
\end{eqnarray}
then
\begin{eqnarray}
|\Delta_{\mu}^{-1}(\gamma_{k+1}^j)\eta^{c_2}|_\alpha\leq
|\Delta_{\mu}^{-1}(\gamma_{k+1}^j)\eta^{c_2}|_{\bar{\alpha}}.
\nonumber
\end{eqnarray}
Thus combining (\ref{eq6.4}) and (\ref{eq6.5}) produces the desired
result (1).

(4) From (1) it immediately follows that
\begin{eqnarray}\label{eeee1}
|w|_\alpha-b\leq |w|_{\bar{\alpha}}\leq |w|_\alpha+b\text{ for
}\alpha\in \Sigma.
\end{eqnarray}
Since $|w|=|w|_\alpha+|w|_{\bar{\alpha}}$, from (\ref{eeee1}) we get
\begin{eqnarray}\label{eeee2}
\frac{|w|}{2}-\frac{b}{2}\leq |w|_\alpha\leq
\frac{|w|}{2}+\frac{b}{2}\text{ for }\alpha\in \Sigma..
\end{eqnarray}
So, combining (\ref{eq4}) and (\ref{eeee2}) gives the desired result
(4).  $\Box$

From Remark 1 and Lemma~\ref{Lem8} (4), we can establish the
following useful bounds of the heights of smooth words of length $n$
for 2-letter even alphabets.
\begin{lem} \label{Lem9} Let $a,b$ be both even numbers. Then  there
are two constants $t_1,t_2$ such that for each positive integer $n$,
ones have
\begin{eqnarray}
ht_{min}(n)&>& \frac{\log n}{\log \rho}+t_1,\label{eeee6}\\
ht_{max}(n)&<& \frac{\log n}{\log \rho}+t_2,\label{eeee5}
\end{eqnarray}
where
\begin{eqnarray}
t_1=-\frac{\log(3b-2+\frac{2(b-1)}{\rho-1})}{\log\rho},\nonumber
\end{eqnarray}
\begin{eqnarray}
t_2=2-\frac{\log (\frac{\rho-2}{\rho-1}b)}{\log \rho}\,\nonumber
\end{eqnarray}
\begin{eqnarray}
\rho=\frac{a+b}{2},\nonumber
\end{eqnarray}
 $-2.3347<-\frac{\log 13}{\log 3}\leq t_1< -1$, $0.7944<2-\frac{\log 20}{\log 12}\leq t_2\leq 2-\frac{\log 2}{\log 3}\approx 1.36907$.
\end{lem}
\noindent{\bf Proof.} First, from  the proof of (\ref{eq2}) and the
right half part of Lemma~\ref{Lem8} (4) it immediately follows the
desired lower bound of $ht_{min}(n)$, where
\begin{eqnarray}
t_1=-\frac{\log(3b-2+\frac{4(b-1)}{a+b-2})}{\log\frac{b+a}{2}}.\nonumber
\end{eqnarray}
Thus
\begin{eqnarray}
t_1< -\frac{\log b}{\log \frac{b+b}{2}}=-1,\nonumber
\end{eqnarray}
and if $a=b-2$ then
\begin{eqnarray}
t_1=-\frac{\log(3b-2+\frac{2(b-1)}{b-2})}{\log (b-1)}\rightarrow -1\,(b\rightarrow \infty).\nonumber
\end{eqnarray}

If $b=4$ then $a=2$, which means $t_1=-\frac{\log 13}{\log 3}$. For
$b\geq 4$, we have
\begin{eqnarray}
t_1&\geq& -\frac{\ln(3b-2+\frac{4(b-1)}{b})}{\ln\frac{b+2}{2}}\nonumber\\
&=&-\frac{\ln\frac{3b^2+2b-4}{b}}{\ln\frac{b+2}{2}}.\label{eqq}
\end{eqnarray}
Let
\begin{eqnarray}
g(b)=\ln 3\ln\frac{3b^2+2b-4}{b}-\ln 13\ln\frac{b+2}{2}, \nonumber
\end{eqnarray}
then
\begin{eqnarray}
g'(b)=(\ln 3)\frac{3b^2+4}{3b^3+2b^2-4b}-\frac{\ln
13}{b+2}.\nonumber
\end{eqnarray}

By Maple, we easily see that the roots of the equation $3(\ln 3-\ln
13)b^3+(6\ln 3-2\ln 13)b^2+4(\ln 3+\ln 13)b+8\ln 3=0$ are
approximately equal to -1.003,\,-0.894,\,2.229. Hence, since
$g'(4)<0$ and $g'(b)$ is continuous in $[4,+\infty)$, we obtain
$g'(b)<0$ for all $b\geq 4$. Therefore $g(b)\leq g(4)=0$, which
suggests
\begin{eqnarray}
\frac{\ln\frac{3b^2+2b-4}{b}}{\ln\frac{b+2}{2}}\leq \frac{\log
13}{\log 3}.\nonumber
\end{eqnarray}
Then (\ref{eqq}) gives $t_1\geq -\frac{\log 13}{\log 3}$.

Second, we use an argument similar to the proof of (\ref{eq2}) to
obtain the upper bound of $ht_{max}(n)$. Note that if $ht(w)\geq 2$,
then $|w|> 2b$. Then from the left half part of Lemma~\ref{Lem8}
(4), we get
\begin{eqnarray}\label{eeee3}
|D(w)|< \frac{1}{\rho} |w|+b.
\end{eqnarray}
Now assume $w$ is a smooth word of length $n$ with height $k$ larger
than or equal to 2. Since $ht(w)\geq 2$, from (\ref{eeee3}), we
arrive at
\begin{eqnarray}
2b&<& |D^{k-2}(w)|\nonumber\\
&<& \frac{1}{\rho}|D^{k-3}(w)|+b\nonumber\\
&<& \frac{1}{\rho^2}|D^{k-4}(w)|+\frac{1}{\rho}b+b\nonumber\\
&\cdots&\nonumber\\
&<&\frac{1}{\rho^{k-2}}|w|+\frac{1}{\rho^{k-3}}b+\cdots
+\frac{1}{\rho^2}b +\frac{1}{\rho}b+b\nonumber\\
&<&\frac{1}{\rho^{k-2}}|w|+\frac{1}{1-\rho^{-1}}b.\nonumber
\end{eqnarray}
Thus
\begin{eqnarray}
\rho^{k-2}<\frac{|w|}{\tau}, \textit{ where
}\tau=\frac{\rho-2}{\rho-1}b,\nonumber
\end{eqnarray}
which means
\begin{eqnarray}
k<\frac{\log n}{\log \rho}+2-\frac{\log \tau}{\log
\rho}.\label{eeee4}
\end{eqnarray}

Note that the length $n$ of a smooth word of height 1 is  greater
than or equal to $a+2\geq 4$, so
\begin{eqnarray}
\frac{\log n}{\log \rho}+2-\frac{\log \tau}{\log \rho}\geq
2+\frac{\log \frac{4(\rho-1)}{b(\rho-2)}}{\log \rho}>1,\nonumber
\end{eqnarray}
which means (\ref{eeee4}) holds for every smooth word. Now from
(\ref{eeee4}) it immediately follows the desired upper bound
(\ref{eeee5}) of $ht_{max}(n)$.

From
\begin{eqnarray}
(b-3)^2-(a-1)^2-8\geq 0\text{ for }b\geq a+4,\nonumber
\end{eqnarray}
we get
\begin{eqnarray}
(b-a)(a+b-4)\geq 2(a+b),\nonumber
\end{eqnarray}
\begin{eqnarray}
2b-(a+b)\geq \frac{2(a+b)}{a+b-4},\nonumber
\end{eqnarray}
\begin{eqnarray}
b\geq \rho+\frac{\rho}{\rho-2},\nonumber
\end{eqnarray}
\begin{eqnarray}
b\frac{\rho-2}{\rho-1}\geq \rho,\nonumber
\end{eqnarray}
which means
\begin{eqnarray}
\frac{\log (\frac{\rho-2}{\rho-1}b)}{\log \rho}> 1>\frac{\log 2
}{\log 3} \text{ for } b\geq a+4.\nonumber
\end{eqnarray}
Thus if $b\geq a+4$ then $t_2\leq 2-\frac{\log 2}{\log 3}\approx
1.36907.$

If $b=a+2$ then $\rho=a+1$, so
\begin{eqnarray}
t_2=2-\frac{\ln \frac{(a-1)(a+2)}{a}}{\ln (a+1)}.\nonumber
\end{eqnarray}
Let
\begin{eqnarray}
f(a)=\ln( 3)\ln \frac{(a-1)(a+2)}{a}-\ln (2)\ln (a+1)\nonumber
\end{eqnarray}
then
\begin{eqnarray}
f'(a)&=&\ln (3)\frac{a^2+2}{(a-1)(a+2)a}-\frac{\ln 2}{a+1}\nonumber\\
&>&\ln (3)(\frac{a^2+2}{(a-1)(a+2)a}-\frac{1}{a+1})\nonumber\\
&=&\ln (3)\frac{4a+2}{(a^2-1)(a+2)a}\nonumber\\
&>&0\text{ for every } a>1.\nonumber
\end{eqnarray}
Hence, $f(a)\geq f(2)=0$ for each $a\geq 2$, that is,
\begin{eqnarray}
\frac{\ln \frac{(a-1)(a+2)}{a}}{\ln (a+1)}\geq \frac{\ln 2}{\ln
3}\text{ for each }a\geq 2,\nonumber
\end{eqnarray}
which also gives the desired result $t_2\leq 2-\frac{\log 2}{\log
3}$.

Finally, machine computation shows
\begin{eqnarray}
t_2\geq 2-\frac{\log 20}{\log 12}\text{ for }b\leq 58.\label{qe1}
\end{eqnarray}
Moreover, in view of $a<b$, we obtain
\begin{eqnarray}
t_2\geq 2-\frac{\log b}{\log \frac{b}{2}}.\label{qe2}
\end{eqnarray}
And let
\begin{eqnarray}
h(b)= 2-\frac{\ln b}{\ln \frac{b}{2}},\label{qe3}
\end{eqnarray}
then
\begin{eqnarray}
h'(b)= \frac{\frac{1}{b}(\ln b-\ln \frac{b}{2})}{(\ln
\frac{b}{2})^2}>0\text{ for }b\geq 4.\nonumber
\end{eqnarray}
which means
\begin{eqnarray}
h(b)\geq h(60)\approx 0.7962>2-\frac{\log 20}{\log 12}\text{ for
}b\geq 60.\nonumber
\end{eqnarray}
Thus (\ref{qe1}), (\ref{qe2}) and (\ref{qe3}) give the desired lower
bound of the constant $t_2$.  $\Box$
\begin{thm} \label{t2} Let $a,b$ be both even numbers. Then  there
exist two suitable constants $c_1,c_2$ such that
\begin{eqnarray}
c_{1} n^{\frac{\log (2b-1)}{\log (a+b)-\log
2}}\leq\gamma_{a,b}(n)\leq c_{2} n^{\frac{\log (2b-1)}{\log
(a+b)-\log 2}},\nonumber
\end{eqnarray}
where $c_1=2(2b-1)^{t_1-1}$, $c_2=2(2b-1)^{t_2}$, $t_1,t_2$ are
determined by Lemma~\ref{Lem9}.
\end{thm}
\noindent{\bf Proof.} From the proof of Theorem~\ref{t1} we easily
see that (\ref{eeq1}) and (\ref{eeq2}) always hold. Thus combining
(\ref{eeee6}) and (\ref{eeq2}) gives
\begin{eqnarray}
\gamma_{a,b}(n) \geq\ c_{1} n^{\frac{\log (2b-1)}{\log (a+b)-\log
2}}.\nonumber
\end{eqnarray}

Similarly, from (\ref{eeee5}) and (\ref{eeq1}) it follows
\begin{eqnarray}
\gamma_{a,b}(n)\leq c_{2} n^{\frac{\log (2b-1)}{\log (a+b)-\log
2}}.\;\;\Box\nonumber
\end{eqnarray}
\section{Concluding remarks\label{s5}}
To establish the estimates of subword complexity function of smooth
words to follow our thoughts and methods is an interesting problem
for large alphabets $\Sigma_n$ containing $n$ letters, where $n\geq
3$.

For the 3-letter alphabet $\Sigma_3=\{2,4,6\}$, let
\begin{eqnarray}
w_1&=&64^22^66^64^66^62^64^6,\nonumber\\
w_2&=&42^66^64^66^62^64^6,\nonumber\\
w_3&=&42^66^64^66^62^6,\nonumber
\end{eqnarray}
\begin{eqnarray}
v_1&=&4^62^26^2,\nonumber\\
v_2&=&2^66^62^66^62^66^64^46^2,\nonumber\\
v_3&=&2^64^62^26^2,\nonumber\\
u_1&=&2^26^24^6,\nonumber\\
u_2&=&4^42^26^22^26^22^26^24^6,\nonumber\\
u_3&=&2^24^26^2,\nonumber
\end{eqnarray}
then $D(w_1)=26^6$, $D(w_2)=6^6$, $D(w_3)=6^5$, $D(v_1)=62$,
$D(v_2)=6^64$, $D(v_3)=6^22$, $D(u_1)=26$, $D(u_2)=2^66$,
$D(u_3)=2$, we easily see that each of $w_1,w_2\text{ and } w_3$ has
only one left smooth extension and $D(w_i)$ has exactly $i$ left
smooth extensions for $i=1,2,3$; each of $v_1,v_2\text{ and } v_3$
has exactly two left smooth extensions and $D(v_i)$ has exactly $i$
left smooth extensions for $i=1,2,3$; each of $u_1,u_2\text{ and }
u_3$ has exactly three left smooth extensions and $D(u_i)$ has
exactly $i$ left smooth extensions for $i=1,2,3$. Thus for large
alphabets containing at least three letters, the estimates of factor
complexity function of smooth words become more complicated than the case for 2-letter alphabets.\\
{\bf Acknowledgements}

I would like to thank Professor William D. Weakley very much for
providing me the opportunity to work with him in \cite{Huang4}.\\


\begin{thebibliography}{WWW}

\bibitem{Ba1} \bf M. Baake, B. Sing, \it Kolakoski-(3,1) is a (deformed) model set, \rm Canad. Math. Bull., 47 (2) (2004), 168-190.

\bibitem{Br1} \bf V. Berth\'{e}, S. Brlek,  P. Choquette, \it Smooth words over arbitrary alphabets,
\rm Theoret. Comput. Sci., 341 (2005), 293-310.

\bibitem{Br2} \bf S. Brlek, S. Dulucq, A. Ladouceur, L. Vuillon, \it Combinatorial properties of
smooth infinite words, \rm Theoret. Comput. Sci., 352 (2006),
306-317.

\bibitem{Br3} \bf S. Brlek, D. Jamet, G. Paquin, \it Smooth
words on 2-letter alphabets having same parity, \rm Theoret. Comput.
Sci., 393 (2008), 166-181.

\bibitem{Br4} \bf S. Brlek, A. Ladouceur, \it A
note on differentiable Palindromes, \rm Theoret. Comput. Sci., 302
(2003), 167-178.

\bibitem{Br5} \bf S. Brlek, G. Melan\c{c}on ,
G. Paquin, \it Properties of the extremal infinite smooth words, \rm
Discrete Math. Theoret. Comput. Sci., 9 (2) (2007), 33-50.

\bibitem{Carpi1} \bf A. Carpi,\, \it Repetitions in
the Kolakovski sequence, \rm Bull. of the EATCS, 50 (1993), 194-196.

\bibitem{Carpi2} \bf A. Carpi,\,
\it On repeated factors in $C^{\infty}$-words, \rm Inf. Process.
Lett., 52 (1994) 289-294.

\bibitem{Chv} \bf V.
Chv\'{a}tal, \it Notes On the Kolakoski sequence, \rm DIMACS Tech.
Rep., 93-84 (1994).

\bibitem{Dekking1} \bf F. M. Dekking,\, \it Regularity
and irragularity  of sequences generated by automata,\, \rm
S\'{e}minaire de Th\'{e}orie des Nombres de Bordeaux, 1979-80,
expos\'{e} n$^{\circ}$ 9, 901-910.

\bibitem{Dekking2} \bf F. M. Dekking,\, \it On the
structure of selfgenerating sequences, \rm S\'{e}minaire de
Th\'{e}orie des Nombres de Bordeaux, 1980-81, expos\'{e} n$^{\circ}$
31, 3101-3106.

\bibitem{Dekking3}
\bf F. M. Dekking,\, \it What is the long range order in the
kolakoski sequence?  \rm in R.V.Moody (ed.), The mathematics of
Long-Range Aperiodic order, Kluwer Academic Publishers (1997)
115-125.

\bibitem{Huang2} \bf Y. B. Huang,\, \it About the number of $C^{\infty}$-words of
form $\tilde{w}xw$, \rm Theoret. Comput. Sci., 393 (2008), 280-286.

\bibitem{Huang3} \bf Y. B. Huang,\, \it The complexity of $C^{b\omega}$-words of the form $\tilde{w}xw$,
\rm Theoret. Comput. Sci., {\bf 410} (2009), 4892-4904.

\bibitem{Huang4} \bf Y. B. Huang, \, W. D. Weakley, \it A note on the complexity of $C^\infty$-words,
\rm Theoret. Comput. Sci., {\bf 411} (2010), 3731-3735.

\bibitem{Huang6} \bf Y. B. Huang,  \it The powers of smooth words over arbitrary 2-letter alphabets,
\rm arXiv:0904.0562v1 (2009).

\bibitem{Keane} \bf M. S. Keane, \,\it Ergodic theory and
subshifts of finite type, \rm in: Ergodic Theory, Symbolic Dynamics
and Hyperbolic Spaces, T. Bedford, M. Keane, C. Series (Eds.),
Oxford University Press (Oxford 1991), 350-370.

\bibitem{Kimberling} \bf C. Kimberling,\, \it Problem 6287, \rm
Amer. Math. Monthly, 86 (1979), 35-70.

\bibitem{Kolakoski} \bf W. Kolakoski,\,\it
Self-generating runs, Problem 5304, \rm Amer. Math. Monthly, 72
(1965), 674. Solution: Amer. Math. Monthly, 73 (1966), 681-682.

\bibitem{Lepist}  \bf A.\, Lepist\"{o},\,\it Repetitions
in the Kolakoski sequence, \rm Development in Language Theory
(1993), 130-143.

\bibitem{P} \bf G. P\v{a}un,\,\it How much Thue is Kolakovski?
\rm  Bull. of the EATCS, 49 (1993), 183-185.

\bibitem{Sing2} \bf B. Sing, \it Kolakoski sequences - an example of aperiodic order, \rm J.
Non-Cryst. Solids, 334-335 (2004), 100-104.

\bibitem{Sing3} \bf B. Sing, \it Kolakoski-(2m,2n) are limit-periodic model sets, \rm J. Math.
Phys. 44 (2) (2003), 899-912.

\bibitem{Steacy} \bf R. Steacy,\,\it Structure in the Kolakoski sequence, \rm
Bull. of the EATCS, 59 (1996), 173-182.

\bibitem{Steinsky} \bf B. Steinsky, \it
A Recursive Formula for the Kolakoski Sequence, \rm  Journal of
Integer Sequences,  9(3) (2006), 06.3.7.

\bibitem{Weakly} \bf W. D. Weakley, \it
On the number of $C^{\infty}$-words of each length, \rm  Jour. of
comb. Theory, Ser.A, 51 (1989), 55-62.

\end{thebibliography}
\end{document}